\documentclass[12pt]{article}

\usepackage{amsmath}
\usepackage{amsthm}
\usepackage{amsfonts}
\usepackage{mathrsfs}
\usepackage{stmaryrd}
\usepackage{setspace}
\usepackage{fullpage}
\usepackage{amssymb}
\usepackage{breqn}
\usepackage{enumitem}
\usepackage{bbold}
\usepackage{authblk}
\usepackage{comment}
\usepackage{hyperref}
\usepackage{pgf,tikz}
\usepackage{graphicx}
\usepackage{subcaption}

\newtheorem{thm}{Theorem}[section]

\newtheorem{lemma}[thm]{Lemma}
\newtheorem{conjecture}[thm]{Conjecture}

\newtheorem{prop}[thm]{Proposition}

\newtheorem{cor}[thm]{Corollary}

\newcommand\ex{\ensuremath{\mathrm{ex}}}

\newcommand\cB{{\mathcal B}}

\newcommand\cF{{\mathcal F}}

\newcommand\cH{{\mathcal H}}

\usepackage{algorithm}      %
\usepackage{algpseudocode}  

\newtheorem{claim}[thm]{Claim}

\newtheorem*{thm*}{Theorem}
\newtheorem*{prop*}{Proposition}
\newcommand{\ignore}[1]{}

\title{On the largest chromatic number of $F$-free hypergraphs}
\author{
Yichen Wang\thanks{\small Department of Mathematical Sciences, Tsinghua University, Beijing, P.R. China. Email:
\small \texttt{wangyich22@mails.tsinghua.edu.cn}.}\,, \hspace{0.2em}
Mengyu Duan\thanks{\small School of Mathematics and Statistics, Northwestern Polytechnical University and Xi'an-Budapest Joint Research Center for Combinatorics, Xi'an 710129, Shaanxi, P.R. China. Email:
\small \texttt{myduan@mail.nwpu.edu.cn}.}\,, \hspace{0.2em}
D\'{a}niel Gerbner\thanks{\small Alfr\'ed R\'enyi Institute of Mathematics. Email:
\small \texttt{gerbner.daniel@renyi.hu}.}\,, \hspace{0.2em} 
Hilal Hama Karim\thanks{\small Department of Computer Science and Information Theory, Faculty of Electrical Engineering and Informatics, Budapest University of Technology and Economics, Műegyetem rkp. 3., H-1111 Budapest, Hungary. E-mail: \texttt{hilal.hamakarim@edu.bme.hu}.}
}

\date{}

\begin{document}

\maketitle

\begin{abstract} 
Given a hypergraph $F$, what is the largest chromatic number that an $F$-free hypergraph can have? In the case of graphs, this question is easy to answer: the chromatic number is unbounded if $F$ contains a cycle, and the largest chromatic number of $F$-free graphs is $k-1$ if $F$ is a forest on $k$ vertices. The situation is more complicated for hypergraphs. 

The strong coloring of a hypergraph is a coloring of the vertices such that every hyperedge is rainbow.
The weak coloring of a hypergraph is a coloring of the vertices such that no hyperedge is monochromatic.
The strong/weak chromatic number of a hypergraph is the minimum number of colors in a strong/weak coloring of the hypergraph. Our question has been completely answered for the weak chromatic number, similarly to the graph case. 

We characterize the hypergraphs $F$ such that $F$-free hypergraphs have bounded strong chromatic number. The only remaining case is when $F$ is the 3-uniform expansion $S_k^+$ of a star with $k$ edges. Concerning the strong chromatic number of $S_k^+$-free hypergraphs, we give bounds that are asymptitically sharp as $k\rightarrow\infty$.

We also consider the same problem when the Berge copies of a graph $F$ are forbidden.
We characterize when the strong/weak chromatic numbers are bounded in this case, and obtain sharp results or bounds for specific trees.
In particular, when $F$ is a path, we give a tight bound when $r=3$ and an asymptotically sharp bound when $r=4$.
\end{abstract}



\section{Introduction}

The study of Tur\'an-type problems forms the cornerstone of extremal graph theory, originating from Tur\'an's generalization of Mantel's theorem \cite{Ma}. What is the largest possible number of edges in $n$-vertex graphs forbidding given subgraphs? The answer is given by Tur\'an~\cite{Tu} in 1941 for complete forbidden graphs. Tur\'an determined the maximum number of edges in a $K_t$-free simple graph and characterized the unique optimal construction, now universally known as the Tur\'an graph. 
In contrast, the Tur\'an problem for $r$-uniform hypergraphs with $r\geq 3$ is complicated, with exact or sharp asymptotic results available only for special classes of forbidden hypergraphs. 
Despite intensive research, it remains one of the most challenging areas of modern extremal combinatorics. 

A hypergraph $\mathcal{H}$ is \textit{$r$-uniform} or simply an \textit{$r$-graph} if it is a 
family of $r$-element subsets of a finite set $V(\mathcal{H})$. 
A hypergraph $\mathcal{H}$ is \textit{at most $r$-uniform} or \textit{$\{2,\ldots,r\}$-uniform} if every hyperedge in $\mathcal{H}$ contains at most $r$ vertices.
Usually, we ignore the hyperedges of size one.
Given a set $\mathcal{F}$ of $r$-graphs, we say a hypergraph is $\mathcal{F}$-free if it contains none of the members of $\mathcal{F}$. We write $\mathrm{ex}_r(n, \mathcal{F})$ for the Tur\'an number of $\mathcal{F}$, i.e.\ the maximum number of edges in an $\mathcal{F}$-free $r$-graph on $n$ vertices. We simply write $\mathrm{ex}_r(n, F)$ if $\mathcal{F} = \{F\}$.

In this paper we consider a related but different extremal parameter: for a fixed hypergraph $F$ we ask how large the chromatic number of an $F$-free hypergraph can be. In the graph case ($r=2$), the answer is simple and well-known.

\begin{thm}[Erd\H{o}s \cite{er}]\label{er}
For any integers $g \geq 3$ and $k \geq 2$, there exists a simple graph $G$ that does not contain cycles of length less than $g$ and has chromatic number at least $k$.
\end{thm}

If $F$ is a forest on $k$ vertices, then we can order its vertices such a way that each vertex has at most one neighbor among the earlier vertices. Therefore, if we try to embed the vertices of $F$ in this order to a graph $G$, we fail only if all the neighbors in $G$ of an already embedded vertex are also already embedded vertices. This cannot happen if every vertex of $G$ has degree at least $k$. Therefore, in an $F$-free graph, we can repeatedly remove a vertex of degree at most $k-1$.
The remaining graph must be empty; otherwise we can get a copy of $F$ greedily.
Then we greedily color the deleted vertices.
Combining Theorem~\ref{er} and the above argument, we get the following result.
\begin{cor}\label{corol}
For a simple graph $F$, if $F$ contains a cycle then there exist $F$-free graphs with arbitrarily large chromatic number; if $F$ is a forest with $k$ vertices then every $F$-free graph has chromatic number at most $k-1$, and $K_{k-1}$ shows this bound is tight.
\end{cor}

    

Note that if the chromatic number is unbounded, that is the final answer for us. However, one could study how large the chromatic number can be as a function of the number of vertices \cite{er}, or the maximum degree \cite{Br}, or the clique number \cite{Gy}, or some other parameter, see \cite{Sc}.

In this paper, we consider the above problem for hypergraphs. In this case, there are multiple notions of the chromatic number. We still consider colorings that partition the vertex set into independent sets (with as few colors as possible), but there are multiple notions of independent sets. The most usual notions are the following. We say that a set $S$ of vertices is \textit{weakly independent} in a hypergraph if no hyperedge is contained in $S$. We say that a set $S$ of vertices is \textit{strongly independent} in a hypergraph if no hyperedge shares two or more vertices with $S$.

Accordingly, weak coloring is such that each color class is weakly independent, while strong coloring is such that each color class is strongly independent.
In other words, in a weak coloring, monochromatic hyperedges are forbidden, while in a strong coloring, each hyperedge is rainbow. 
The \textit{strong chromatic number}~(resp. \textit{weak chromatic number}) of $\mathcal{H}$, denoted by $\chi_s(\mathcal{H})$~(resp. $\chi_w(\mathcal{H})$), is the minimum number of colors such that it admits a strong coloring~(resp. weak coloring).
Note that when $\mathcal{H}$ is $r$-uniform, if we unite disjoint $(r-1)$-sets of color classes in a strong coloring, then we obtain a weak coloring, thus $\chi_w(\mathcal{H})\le \lceil\chi_s(\mathcal{H})/(r-1)\rceil$.

For a set $\mathcal{F}$ of hypergraphs, we define
\[
\mathrm{schex}_r(n, \mathcal{F})=\sup\{\chi_s(\mathcal{H}): \mathcal{H} \text{ is an $\mathcal{F}$-free $r$-uniform hypergraph on $n$ vertices}\},
\]
\[
\mathrm{wchex}_r(n, \mathcal{F})=\sup\{\chi_w(\mathcal{H}): \mathcal{H} \text{ is an $\mathcal{F}$-free $r$-uniform hypergraph on $n$ vertices}\},\]and then it can be $\infty$ when $n$ goes to infinity.

We also investigate this problem for Berge hypergraphs. 
Let $F$ be a graph. We say that a hypergraph $H$ is a \emph{Berge--$F$}
if there is a bijection
$f : E(F) \to E(H)$ such that $e \subseteq f(e)$ for every $e \in E(F)$. We call the vertices and edges of $F$ and the hyperedges of $H$ \textit{defining} vertices, edges and hyperedges.
Note that there are multiple non-isomorphic hypergraphs that are Berge--$F$, let $\mathcal{B}(F)$ denote their family.

The concept of Berge hypergraphs in the case of cycles was first introduced by Claude Berge in his seminal 1973 monograph on hypergraph theory \cite{Be}, as a way to generalize graph concepts to the hypergraph setting by preserving the incidence relation between vertices and edges. Gerbner and Palmer \cite{B3} generalized it to arbitrary graphs $F$. The Berge--$F$ construction, in particular, defines a broad family of hypergraphs that inherit the structural essence of a fixed base graph $F$: each edge of $F$ is mapped to a hyperedge containing it, via a bijection between the edge sets. This framework has become a central object of study in extremal hypergraph theory, with research dedicated to Tur\'an-type problems for Berge--$F$-free hypergraphs, including tight bounds on their maximum edge count for numerous classes of base graphs $F$. In this work, we extend this line of inquiry by analyzing the extremal chromatic number of Berge--$F$-free hypergraphs, complementing existing results on edge-maximizing extremal constructions. Note that the strong chromatic number of a hypergraph is equal to the chromatic number of its shadow graph (the graph consisting of the edges contained in some hyperedges). Therefore, our results fit into an existing line of research \cite{luwan}, where we study extremal properties of the shadow graph of Berge-$F$-free hypergraphs.

We will consider two specific Berge hypergraphs, mostly to simplify notations.  Let $F$ be a graph. The \emph{$r$-expansion} of $F$, denoted by $F^{(r)+}$, is obtained from $F$ by enlarging each edge of $F$ with a set of $r-2$ new vertices such that all these new sets are disjoint; the \emph{$r$-suspension} of $F$, denoted by $S^{r}(F)$, is the $r$-graph obtained from $F$ by taking $r-2$ new vertices and adding them to every edge. 

The weak chromatic number of single hypergraphs is already well-studied. Let us call a hypergraph \textit{hypertree} if it does not contain any Berge cycle (not even a Berge cycle of length two, i.e., two hyperedges sharing two vertices).
Loh~\cite{loh}, strengthening a result of Bohman, Frieze and Mubayi~\cite{bfm} proved that every $r$-uniform hypergraph with weak chromatic number larger than $k$ contains a copy of every $r$-uniform hypertree with $k$ edges.
In our language, it shows that $\mathrm{wchex}_r(n, T) \le k$ for any $r$-uniform hypertree $T$ with $k$ edges. Note that this bound is obviously sharp, since $T$ has $(r-1)k+1$ vertices, thus the complete $r$-uniform hypergraph on $(r-1)k$ vertices does not contain $T$ and has weak chromatic number $k$.
Gy\'{a}rf\'{a}s and Lehel~\cite{GYARFAS2011208} pointed out that the greedy algorithm for graphs we described before Corollary \ref{corol} also gives the upper bound $k$ in this case. 
Therefore, in this paper, we mainly focus on the strong chromatic number.

We have the following characterization of $\mathrm{schex}_r(n, \mathcal{F})$ and $\mathrm{wchex}_r(n, \mathcal{F})$ when $\mathcal{F}$ is a single hypergraph or a Berge hypergraph family.

\begin{prop}\label{propp}
    \textbf{(i)} Let $r\ge 3$ and $F$ be an $r$-graph. Then the strong chromatic number of $F$-free $r$-graphs is bounded if and only if either $F$ has at most one edge, or $r=3$ and $F$ is the expansion of a star.

    \textbf{(ii)} Let $r\ge 3$ and $F$ be an $r$-graph. Then the weak chromatic number of $F$-free $r$-graphs is bounded if and only if $F$ does not contain any Berge cycle.

    \textbf{(iii)} Let $r\ge 3$ and $F$ be a graph. Then the strong chromatic number and the weak chromatic number of Berge-$F$-free $r$-graphs are bounded if and only if $F$ does not contain any Berge cycle.
\end{prop}

Note here again that we stop when we get unbounded chromatic number.
But it could be studied how large the chromatic number can be as a function of other parameters.

For $\mathrm{schex}_r(n, \mathcal{F})$ when $\mathcal{F}$ only contains a single hypergraph, it is bounded and non-trivial only when $r=3$ and $\mathcal{F} = \{F\}$ is the expansion of a star. 
In this case, we have the following theorem.

\begin{thm}\label{thm: main single expansion star r=3}
We have $\mathrm{schex}_3(n, S_k^+)=(1+o(1))2k^2$ as $k$ goes to infinity.
  More precisely, we have $\mathrm{schex}_3(n, S_k^+)\le 2k^2+k/2-1$. Moreover, for infinitely many $k$, we have $\mathrm{schex}_3(n, S_k^+)\ge 2k^2-2k-2$.
\end{thm}

We then focus on the case for Berge--$F$-free hypergraphs when $F$ is a forest.
First we have the following general upper bound.
\begin{thm}\label{thm: main Berge-T, general upper bound}
  Let $T$ be a forest with $k+1$ vertices and maximum degree $\Delta$.
  Then we have the following:
  \begin{equation*}
  \begin{aligned}
    \mathrm{schex}_r(n, \mathcal{B}(T)) &\le k+(r-3)(k-1)+\Delta -1.
  \end{aligned}
  \end{equation*}
\end{thm}

When $T = S_k$ is a star with $k$ edges, 
Theorem~\ref{thm: main Berge-T, general upper bound} gives $\mathrm{schex}_r(n, \mathcal{B}(S_k)) \le (r-1)(k-1)+1$.
The Steiner system $\mathcal{S}(n, r, t)$ is an $n$-vertex $r$-uniform hypergraph such that every $t$-set is covered exactly once.
When the Steiner system $\mathcal{S}((r-1)(k-1)+1, r, 2)$ exists, denote it by $\mathcal{H}$.
Consider the strong coloring of $\mathcal{H}$.
Since every pair is contained in a hyperedge, every pair receives different colors, which means $\chi_s(\mathcal{H}) \ge (r-1)(k-1)+1$.
It is easy to see that $\mathcal{H}$ is $\mathcal{B}(S_k)$-free.
Then the upper bound from Theorem~\ref{thm: main Berge-T, general upper bound} in the strong coloring case is tight for $T=S_k$.

When $T = P_k$ is a path with $k$ edges and $r=3$, we have a better characterization.
\begin{thm}\label{thm: main Berge path r=3}
  When $k \ge 3$, we have 
  \begin{equation*}
  \begin{aligned}
    \mathrm{schex}_3(n, \mathcal{B}(P_k)) &= k, \\
    \mathrm{wchex}_3(n, \mathcal{B}(P_k)) &= \left\lceil \frac{k}{2} \right\rceil.
  \end{aligned}
  \end{equation*}
\end{thm}

Moreover, we will present an algorithm to find a strong coloring of $\mathcal{B}(P_k)$-free $3$-graphs with $k$ colors and a weak coloring of $\mathcal{B}(P_k)$-free $3$-graphs with $\lceil \frac{k}{2} \rceil$ colors~(see Theorem~\ref{thm: DFS algorithm strong r=3 path} and Theorem~\ref{thm: weak coloring DFS algorithm r=3 path}).

We conjecture that 
\begin{conjecture}\label{conj: main berge path}
  Let $r \ge 3$ be a fixed integer.
  There exists a function $f(r)$ such that when $k \ge f(r)$, we have 
  \begin{equation*}
    \mathrm{schex}_r(n, \mathcal{B}(P_k)) = k.
  \end{equation*}
\end{conjecture}

Note that a large enough $k$ is a necessary condition.
In the $4$-uniform case, we prove a weaker result:
\begin{thm}\label{thm: main 4uniform berge path}
  \begin{equation*}
    \mathrm{schex}_4(n, \mathcal{B}(P_k)) \le k + 11.
  \end{equation*}
\end{thm}

In Subsection~\ref{subsec: 3uniform special trees}, we determine $\mathrm{schex}_3(n, \mathcal{B}(T))$ for some special trees $T$ like spider, double star and broom~(see Theorem~\ref{thm: spider}, Theorem~\ref{thm: double star} and Theorem~\ref{thm: broom}).

The paper is organized as follows.
In Section~\ref{sec: preliminary}, we introduce some basic definitions and notations.
In Section~\ref{sec: proof of propp}, we prove Proposition~\ref{propp} and Theorem~\ref{thm: main single expansion star r=3}.
In Section~\ref{sec: berge}, we study the strong/weak chromatic number for Berge-tree-free hypergraphs in the $3$-uniform case.
In Section~\ref{sec: berge path}, we prove Theorem~\ref{thm: main Berge path r=3} and Theorem~\ref{thm: main 4uniform berge path}.

\section{Preliminaries}\label{sec: preliminary}

A hypergraph is \emph{linear} if every pair of hyperedges intersect in at most one vertex. 
The \emph{2-shadow} or \emph{shadow graph} of a hypergraph $\cH$, denoted by $\partial \cH$, is the graph on $V(\cH)$ where $uv\in E(\partial \cH)$ if and only if there exists $e\in E(\cH)$ with $\{u,v\}\subseteq e$. $N_\cH(v)$ denotes the set of vertices that are contained in a hyperedge together with $v$, i.e., the \emph{neighbors} of $v$.
The \emph{link hypergraph} of $v$ in $\mathcal{H}$ denoted by $L_{\mathcal{H}}(v)$, is the hypergraph on vertex set $N_{\cH}(v)$ where $e$ is a hyperedge if and only if $e\cup \{v\}$ is a hyperedge in $\mathcal{H}$. 
Especially, when $\mathcal{H}$ is $r$-uniform, the link hypergraph of $v$ is $(r-1)$-uniform.

To extend the simple argument in the graph case (described before Corollary \ref{corol}) to hypergraphs, we consider three kinds of degrees. 
Given a hypergraph $\cH$ and a vertex $v$ of $\cH$, $d(v)=d_\cH(v)$ denotes the number of hyperedges in $\cH$ that contain $v$. 
We denote by $d^N(v)=d^N_\cH(v)$ the number of vertices that are contained in a hyperedge of $\cH$ together with $v$, that is, $d^N(v)=\left| N_{\mathcal{H}}(v) \right|$. 
We denote by $d^M(v)=d^M_\cH(v)$ the size of the largest matching in the link graph of $v$.
Here a matching is a set of disjoint edges.

For a hypergraph $\mathcal{H}$ and a vertex $u$ in $\mathcal{H}$, let $\mathcal{H}\setminus u$ be the hypergraph obtained from $\mathcal{H}$ by removing $u$ from every hyperedge containing $u$ (without changing the other hyperedges), that is, a hypergraph with vertex set $V(\mathcal{H})\setminus \{u\}$ and hyperedge set $\{e \setminus \{u\} \mid e \in E(\mathcal{H})\}$. 

For a tree $T$ rooted at $r$,
the \textit{depth} of a vertex $v$ is the number of edges on the path from $r$ to $v$. 
Especially, the depth of $r$ is zero.
The \textit{height} of $T$ is the maximum depth of a vertex in $T$.

In the proof, very often we need to find a Berge copy of $T$ in a hypergraph.
For a bijection $f: E(T) \to E(\mathcal{H})$ such that $e \subseteq f(e)$ for every $e \in E(T)$, we say we \textit{embed} $T$ in $\mathcal{H}$ using $f$ on the vertex set of $T$.
Then for a vertex $v$ in $T$, it is also a vertex in $\mathcal{H}$. 
The vertices of $T$ are called the \textit{defining vertices} and every $f(e)$ is called a \textit{defining hyperedge}.




\section{Proof of Proposition~\ref{propp} and Theorem~\ref{thm: main single expansion star r=3}}\label{sec: proof of propp}

We will use the following result. 

\begin{thm}[Erd\H{o}s and Hajnal~\cite{eha}]\label{erha}
For any integers $g \geq 3$, $k \geq 2$, and $r \geq 3$, there exists an $r$-uniform hypergraph $H$ that does not contain Berge cycles of length less than $g$ and has weak chromatic number at least $k$.
\end{thm}

Now we are ready to prove Proposition~\ref{propp}. Recall that it states that $\mathrm{schex}_r(n, F)$ is bounded if and only if $F$ contains only one edge or $r=3$ and $F=S_k^+$ for some $k$, while $\mathrm{schex}_r(n, \mathcal{B}(F))$ is bounded if and only if $\mathrm{wchex}_r(n, \mathcal{B}(F))$ is bounded, which holds if and only if $F$ is a forest.

\begin{proof}[Proof of Proposition \ref{propp}] 
  Let us start with \textbf{(i)} and assume first that $F$ has at least two edges, and let $t$ be the size of the intersection of two edges. If $t<r-2$, then the $r$-suspension of $K_k$ has strong chromatic number at least $k$, and is $F$-free since any two hyperedges share at least $r-2$ vertices. If $t>1$, then the $r$-expansion  of $K_k$ has strong chromatic number at least $k$, and is $F$-free since any two hyperedges share at most 1 vertex. This proves that $\mathrm{schex}_r(n,F)$ is unbounded for any $F$ unless $r=3$ and $F$ is linear. If $F$ is not the expansion of a star, consider the suspension of $K_k$. Since each hyperedge contains the suspending vertex $v$, the only linear hypergraphs contained in this hypergraph contain $v$ and its hyperedges are disjoint otherwise. But this means that our hypergraph is the expansion of a star. 
  
  If $F$ has at most one hyperedge, then an $F$-free $r$-graph either has fewer vertices than $F$, or no hyperedges. The chromatic number is clearly bounded in both cases. Finally, if $r=3$ and $F$ is the expansion of a star, then Theorem~\ref{thm: main single expansion star r=3} completes the proof.

      Let us continue with \textbf{(ii)} and \textbf{(iii)} together. Let us assume that $F$ contains a Berge cycle of length $k$. Then for every $m$, Theorem \ref{erha} provides an $F$-free and Berge-$F$-free $r$-graph with weak chromatic number at least $m$, thus also strong chromatic number at least $m$.

  If $F$ is an $r$-uniform hypertree, then the result of Loh~\cite{loh} completes the proof.
      Finally, consider Berge copies of a forest $F$ on $k$ vertices. 
      We will use a similar argument. We claim that any Berge-$F$-free $r$-graph contains a vertex $u$ with $d(u)\le \binom{k-2}{r-1}+k-1$. Indeed, first we extend $F$ to a $k$-vertex tree $F'$ by adding edges that connect the components, arbitrarily. Recall that any tree can be built from a single edge by adding edges one by one, using the following steps: we take a vertex $v$ already in the tree, and add an edge joining $v$ to a new vertex. 
      If we are given an $r$-graph where each vertex has degree at least $\binom{k-2}{r-1}+k$, then we embed the first edge of $F'$ to two vertices of an arbitrary hyperedge, and then add the further edges in the above described order. At each step, we are given an already embedded defining vertex $u$, and we want to embed an edge incident to $u$. There are at most $k-1$ already embedded hyperedges, those cannot be used. Furthermore, we cannot use a hyperedge if it contains only defining vertices. There are at most $\binom{k-2}{r-1}$ such hyperedges, thus we can find the desired hyperedge. This way we embed $F'$, thus also $F$.
  
      We claim that the strong chromatic number of Berge-$F$-free $r$-graphs is at most $(r-1)\left(\binom{k-2}{r-1}+k-1\right)+1$.  We apply induction on $n$, the case $n\le (r-1)\left(\binom{k-2}{r-1}+k-1\right)+1$ is trivial. For arbitrary $n$, we take a vertex $u$ such that $d(u)\le (r-1)\left(\binom{k-2}{r-1}+k-1\right)$, remove $u$, and color the rest. After that, a color is forbidden for $u$ if and only if there is a vertex of that color in the neighborhood of $u$. Hence, we find a color for $u$ that is not forbidden. 
      
      This implies the same upper bound on the weak chromatic number of Berge-$F$-free $r$-graphs and completes the proof.
  \end{proof}
  

  Now let us prove Theorem~\ref{thm: main single expansion star r=3}.
  We will use the following result of  Gy\'{a}rf\'{a}s,  Lehel, Ne{\v{s}}et{\v{r}}il, R{\"o}dl, Schelp, and Tuza~\cite{gylnrst}.

  \begin{thm}[Gy\'{a}rf\'{a}s,  Lehel, Ne{\v{s}}et{\v{r}}il, R{\"o}dl, Schelp, and Tuza \cite{gylnrst}]\label{locaco}
      Assume that in an edge-colored graph $G$, there is no rainbow $S_k$. Then at least one of the colors has at least $|E(G)|/(k-1)$ edges.
  \end{thm}

  We also use the following result on the Tur\'an number of matchings.

  \begin{thm}[Erd\H{o}s and Gallai \cite{EGa}]\label{thm: erdos gallai}
     If $n\ge 2k-1$, then
      $\ex(n,M_k)=\max\{\binom{2k-1}{2}, \binom{k-1}{2}+(k-1)(n-k+1)\}$.
  \end{thm}

  Now we are ready to prove Theorem~\ref{thm: main single expansion star r=3} that we restate here for convenience.

\begin{thm*}
We have $\mathrm{schex}_3(n, S_k^+)=(1+o(1))2k^2$ as $k$ goes to infinity.
  More precisely, we have $\mathrm{schex}_3(n, S_k^+)\le 2k^2+k/2-1$. Moreover, for infinitely many $k$, we have $\mathrm{schex}_3(n, S_k^+)\ge 2k^2-2k-2$.
\end{thm*}

  \begin{proof}
    Let us start with the upper bound.
    Let $\cF$ be a counterexample with minimum size and $c=2k^2+k/2-1$.
      Denote by $G$ the shadow graph of $\cF$. 
      If $G$ has a vertex $u$ of degree at most $c-1$, then we delete such a vertex. 
      We repeat this as long as we can, thus ending with a subgraph $G'$ with minimum degree at least $c$. 
      Let $U=V(G')$, $m=|U|$ and let $\cF'$ denote the hypergraph consisting of the hyperedges inside $U$. 
      Note that $G'$ must be non-empty, otherwise we can color the deleted vertices greedily. 
      Thus, $m > c$ by the definition.
      
      Recall that $\cF'$ is $S_k^+$-free. Therefore, for each vertex $u\in U$, its link graph $L_{\cF'}(u)$ in $\cF'$ is $M_k$-free. Let us denote by $S(u)$ a smallest covering set in $L_{\cF'}(u)$. Clearly, $|S(u)|\le 2k-2$, for example, the vertices of a largest matching form a covering set. We delete from $G'$ all the edges that join a vertex $u$ to a vertex in $S(u)$, and denote by $G''$ the resulting graph. 
  
      Consider a hyperedge $uvw$ in $\cF'$. We claim that at most one of its subedges is in $G''$. Assume that $uv$ and $uw$ are in $G''$. Then $vw$ is an edge of $L_{\cF'}(u)$, thus at least one of $v$ and $w$ is in $S(u)$, hence we deleted at least one of the edges $uv$, $uw$, a contradiction.
  
      Since $G'$ has at least $cm/2$ edges, $G''$ has at least $cm/2-(2k-2)m$ edges. Each edge $uv$ of $G''$ is contained in a hyperedge $uvw$ of $\cF$; we say that the color of $uv$ is $w$ (we pick $w$ arbitrarily in the case $uv$ is contained in more than one hyperedge of $\cF$). This way, we colored the edges with at most $n$ colors. 
      
      We claim that there is no rainbow $S_k$, nor monochromatic $M_k$ in $G''$. Assume that there is a rainbow $S_k$ with center $u$ and leaves $v_1,\dots, v_k$ in $G''$. Then for each edge $uv_i$ there is a vertex $w_i$ such that $uv_iw_i$ is a hyperedge in $F$, and the vertices $w_i$ are pairwise distinct. These hyperedges do not form $S_k^+$, hence we have $w_i=v_j$ for some $j$. Then for the hyperedge $u_iv_iv_j$ in $F'$, both $uv_i$ and $uv_j$ are in $G''$, a contradiction. A monochromatic $M_k$ clearly means an $M_k$ in the link graph of a vertex, again a contradiction.

      By Theorem~\ref{locaco}, at least one of the colors has at least $(cm/2-(2k-2)m)/(k-1)$ edges. On the other hand, the number of edges in any color is at most $\max\{\binom{2k-1}{2}, \binom{k-1}{2}+(k-1)(m-k+1)\}$ by the Erd\H{o}s--Gallai theorem~(Theorem~\ref{thm: erdos gallai}). Since $m>c$, it is easy to see that $\binom{k-1}{2}+(k-1)(m-k+1)$ is larger, thus we have $\binom{k-1}{2}+(k-1)(m-k+1)\le (cm/2-(2k-2)m)/(k-1)$, a contradiction.

  Let us continue with a construction for the lower bound. 
  Assume that $k-1$ is a prime power. We take a finite projective plane of order $k-1$, with lines $\ell_1,\dots, \ell_{k^2-k-1}$ and vertices $v_1, \dots, v_{k^2-k-1}$. Then we take another copy, with lines $\ell_1',\dots, \ell_{k^2-k-1}'$ and vertices $v_1', \dots, v_{k^2-k-1}'$. Let $b_i=\ell_i\cup \ell_i'$. We add additional vertices $u_1, \dots, u_{k^2-k-1}$ and take as hyperedges the 3-sets that contain some $u_i$ and two vertices of $b_i$ for some $i$. 
  
  A vertex $v_i$ is contained in $k-1$ sets $b_j$. Thus each hyperedge that contains $v_i$ also contains one of the $k-1$ corresponding vertices $u_j$. Therefore, $v_i$ is not the center of an $S_k^+$, and the same holds for $v_i'$. For a vertex $u_i$, there are $2k-2$ vertices that are in a hyperedge together with $u_i$, hence $u_i$ is also not a center of an $S_k^+$. This shows that the resulting hypergraph is $S_k^+$-free.

  We claim that the chromatic number of this hypergraph is at least $2k^2-2k-2$. For $v_i$ and $v_j$ there is a line $\ell_p$ containing both $v_i$ and $v_j$. Therefore, there is a hyperedge $v_iv_ju_p$ containing both of these vertices. Similarly, the pairs $v_i$ and $v_j'$, $v_i'$ and $v_j$, $v_i'$ and $v_j'$, and also $v_i$ and $v_i'$ each form a hyperedge with $u_p$. Therefore, all the $2k^2-2k-2$ vertices $v_1, \dots, v_{k^2-k-1}, v_1', \dots, v_{k^2-k-1}'$ have different colors.

Finally, if $k$ is not a prime power, let $k_0$ be the largest prime power below $k$; the same construction (with $k_0$ in place of $k$) shows that $\mathrm{schex}_r(n, S_k^+)\ge 2k_0^2-2k_0-2$. We obtain the asymptotic lower bound by observing that $k_0$ is sufficiently close to $k$, for example, using a result of Baker, Harman and Pintz \cite{bhp} that there is a prime between $x$ and $x-x^{0.525}$.

  \end{proof}

\section{Strong chromatic number for Berge-tree-free hypergraphs}\label{sec: berge}

In this section, we study the strong and weak chromatic number of Berge-tree-free hypergraphs and prove Theorem~\ref{thm: main Berge-T, general upper bound}.
For some special trees, like spider, double star and broom, we can get a sharp bound~(see Theorem~\ref{thm: spider}, Theorem~\ref{thm: double star} and Theorem~\ref{thm: broom} respectively).

\subsection{Berge-tree-free hypergraphs considering maximum degree}\label{subsec: r uniform tree}

First we prove Theorem~\ref{thm: main Berge-T, general upper bound} in this subsection. 
We actually prove a stronger result which implies Theorem~\ref{thm: main Berge-T, general upper bound} as follows:

\begin{prop}\label{prop: berge tree, r uniform}
  Let $\mathcal{H}$ be an at most $r$-uniform hypergraph without Berge copies of $T$ where $T$ is a tree with $k$ edges and maximum degree $\Delta$.
  Then $\chi_s(\mathcal{H}) \le  k+(r-3)(k-1)+\Delta-1$.
  \end{prop}

  The following lemma contains the key greedy idea used in the whole paper.

\begin{lemma}\label{lem: minimum degree to berge T}
    Let $T$ be a tree with $k$ edges, let $\Delta$ be the maximum degree of $T$.
    Let $\mathcal{H}$ be an at most $r$-uniform hypergraph.
    If for every vertex $u$, we have $d^N(u) \ge k+(r-3)(k-1)+\Delta-1$, then $\mathcal{H}$ contains a Berge copy of $T$.
\end{lemma}

\begin{proof}
    We would like to embed the Berge copy of $T$ step by step, in the order such that each vertex has at most one earlier neighbor, just like in the graph case described before Corollary \ref{corol}.
    Suppose some vertices in $T$ are already embedded and now we are trying to embed the vertex $v$ such that there is an edge $uv$ in $T$ and $u$ is already embedded to a vertex $u'$ of $\cH$.
    Since $d^N(u') \ge k+(r-3)(k-1)+\Delta-1$, we have this many candidate vertices for embedding $v$.
    Let $u_1,u_2,\ldots,u_m$ be the neighbors of $u$ which are embedded.
    Since $v$ is not yet embedded and the maximum degree of $T$ is $\Delta$, we have $m \le \Delta-1$.
    Then we claim that there exists a vertex $v' \in N_{\mathcal{H}}(u')$ such that 
    \begin{enumerate}
      \item no vertex is embedded to $v'$,
      \item for every defining hyperedge $e'$ of an embedded edge such that $u'\in e'$, $v'$ is not in $e'$.
    \end{enumerate}

    Let $A$ be the set of embedded vertices in $T$, then $|A| \le k$. 
    Let $B_1$ be the set of vertices contained in a defining hyperedge $e'$ of an embedded edge $e$ such that $u \in e$, then $|B_1\setminus A| \le (r-2)m$.
    Let $B_2$ be the set of vertices contained in a defining hyperedge $e'$ of an embedded edge $e$ such that $u \notin e$, but $u' \in e'$, then $|B_2\setminus A| \le (r-3)(k-1-m)$. Indeed, we pick the at most $k-1-m$ defining hyperedges that do not correspond to an edge $uu_i$, and if they also contain $u'$, then they contain at most $r-3$ vertices not in $A$.
    
    Therefore, we can prove the existence of $v'$ with the prescribed properties. 
    Then we can embed $v$ to $v'$ such that $u'v'$ is in a hyperedge that is not a defining hyperedge.
    
    We repeat the process and eventually we have embedded a Berge copy of $T$ in $\mathcal{H}$.
\end{proof}

The above Lemma easily implies Proposition~\ref{prop: berge tree, r uniform}.

\begin{proof}[Proof of Proposition~\ref{prop: berge tree, r uniform}]
  We prove the proposition by induction on the number of vertices in $\mathcal{H}$.
  By Lemma~\ref{lem: minimum degree to berge T}, when $d^N(u) \ge k+(r-3)(k-1)+\Delta-1$ for every vertex $u$, we are done.
  So we may assume that there exists a vertex $u$ such that $d^N(u) < k+(r-3)(k-1)+\Delta-1$.

  Then let $\mathcal{H}' = \mathcal{H}\setminus \{u\}$.
  By the induction hypothesis, $\mathcal{H}'$ admits a strong coloring with $k+(r-3)(k-1)+\Delta-1$ colors.
  Then we can greedily color $u$ with a color that is not used in $N_{\partial \mathcal{H}}(u)$.
\end{proof}

As mentioned in the introduction, when $T=S_k$, the bound is tight.

\subsection{Forbidding specific Berge trees: spider, double star and broom }\label{subsec: 3uniform special trees}

In this subsection, we determine $\mathrm{schex}_3(n, \mathcal{B}(T))$ for some special trees $T$.

\begin{figure}[t]
  \centering
  \includegraphics[width=\textwidth]{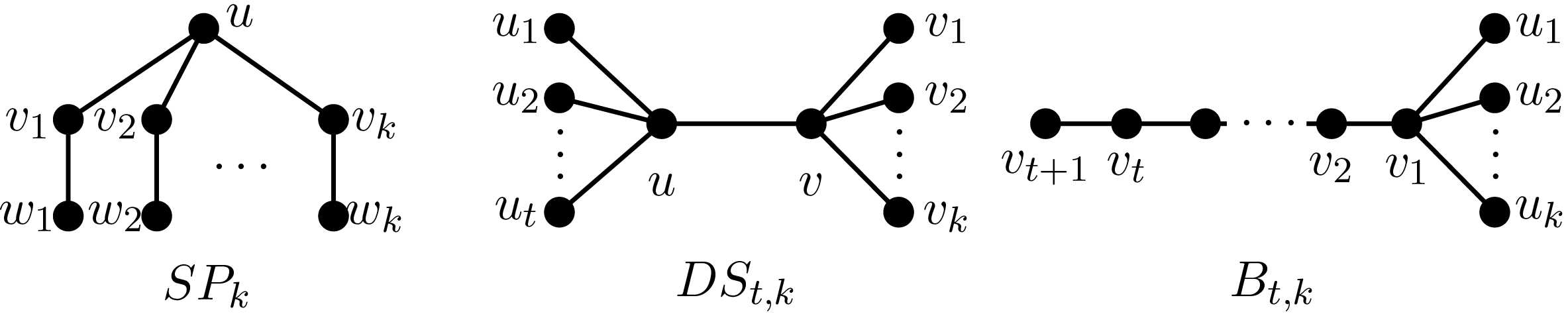}
  \caption{The spider $SP_k$, the double star $DS_{t,k}$ and the broom $BR_{t,k}$.}\label{fig: special trees}
\end{figure}

\begin{thm}\label{thm: spider}
Let $T$ be a spider with $k$ legs of length two~(see Figure~\ref{fig: special trees}), denoted by $SP_k$, then 
\[
  \mathrm{schex}_{3}(n, \mathcal{B}(SP_k)) = 2k.
\]
\end{thm}

\begin{proof}
  The lower bound is obvious by the construction of a clique on $2k$ vertices.
  Let us focus on the upper bound.
  Like in Section~\ref{subsec: r uniform tree}, we actually prove a stronger version for at most $3$-uniform hypergraphs.
  That is, for every at most $3$-uniform hypergraph $\mathcal{H}$ without any Berge copy of $SP_k$, we have $\chi_s(\mathcal{H}) \le 2k$.
  Let $\mathcal{H}$ be an at most $3$-uniform counter-example with the minimum number of hyperedges. 
  
  We may assume $d^N(u) \ge 2k$ for every vertex $u$ of $\cH$.
  Otherwise, let $v$ be a vertex with $d^N(v) < 2k$, then there is a coloring of $\mathcal{H}\setminus \{v\}$ with $2k$ colors by our assumption on $\cH$ and then we can color $v$ greedily.

  Now we would like to find a Berge copy of $SP_k$ in $\mathcal{H}$. We will embed the vertices of $SP_k$ one by one, and each time when we complete an edge $xy$ by embedding its second vertex, we will pick the corresponding hyperedge $f(xy)$.
    At first, we map the center vertex $u$ to an arbitrary vertex $u'$.
    Since $d^N(u') \ge 2k$, there exist vertices $v'_1,v'_2\ldots,v'_k$ such that each $u'v_i'$ is contained in a distinct hyperedge. We map $v_i$ to $v_i'$ and pick the corresponding hyperedge as $f(u'v'_i)$.

    Then we try to embed the edges $v_1w_1, v_2w_2,\ldots,v_{k-1}w_{k-1}$. 
  By a similar argument as in the proof of Proposition~\ref{prop: berge tree, r uniform}, there exists $w'_1,w'_2,\ldots,w'_{k-1}$ such that each edge $v'_iw'_i,i=1,2,\ldots,k-1$ is contained in a distinct hyperedge. Indeed, if $v'_i$ is contained in an already embedded hyperedge $e$ that is not $f(u'v'_i)$, then the other two vertices of $e$ are already embedded. Since we have embedded at most $2k-2$ vertices other than $v_i$, there are at least two unused neighbors of $v'_i$. One of those may be contained in $f(u'v'_i)$, but the other is contained in an unused hyperedge together with $v'_i$.

    Then it remains to embed the edge $v_kw_k$.
    The same approach only fails when 
    \[
    N_{\mathcal{H}}(v'_k) = \{u', v'_1,\ldots,v'_{k-1},w'_1,\ldots,w'_{k-1}, s_{k}\},
    \]
    where $s_k$ is the third vertex in the defining hyperedge of $uv_k$.

    We restart this procedure with an order where $v_i$ is the last vertex. Then we obtain a vertex $s_i$ similarly. We do this for every $i\le k$. If $s_i=s_j$ for some $i\neq j$, then we map $v_j$ to $s_j$ instead of $v_j'$, and we still pick the hyperedge $u'v_j's_j$ as the corresponding defining hyperedge. Then we continue the above way in the order where $v_i$ is the last vertex. Then $N_{\mathcal{H}}(v'_i) = \{u', v'_1,\ldots,v'_{i-1}, v'_{i+1},\ldots, v_k,w'_1,\ldots,w'_{i-1}, w'_{i+1},\ldots, w_k, s_{i}\}$, but then $v_i'$ has only $2k-1$ neighbors, a contradiction. 
    
    Therefore, $s_i \neq s_j$ for $i \neq j$ and $\{s_1,\ldots,s_k\}$ and $\{v'_1,\ldots,v'_k\}$ are disjoint.
 Moreover, for each $i$, $s_i$ and $v'_i$ are equivalent, so $\{u,v'_1,\ldots,v'_k,s_1,\ldots,s_k\}$ is a $(2k+1)$-clique in $\partial \mathcal{H}$.
    Then we map $w_i$ to $ s_{i+1}$ for $ i=1,2,\ldots,k$ where $s_{k+1}=s_1$.
    Each $v'_is_{i+1}$ is contained in a hyperedge since it is an edge in $\partial \mathcal{H}$.
    The hyperedge is not the defining hyperedge of any $u'v'_i$ and for each $v'_is_{i+1}$, the hyperedge is different.
    Then we have a Berge copy of $SP_k$, a contradiction.
\end{proof}

\begin{thm}\label{thm: double star}
Let $T$ be a double star $DS_{t, k}$ with $t+k+2$ vertices, that is, an $S_{t+1}$ and an $S_{k+1}$ sharing one edge~(see Figure~\ref{fig: special trees}). When $t \le k$,
\[
    \mathrm{schex}_3(n, \mathcal{B}(DS_{t,k})) = 2k+1.
\]
\end{thm}

\begin{proof}
  Since $S_{k+1}$ is a subgraph of $DS_{t,k}$, we have 
  \[
        2k+1 \le \mathrm{schex}_3(n, \mathcal{B}(S_{k+1})) \le \mathrm{schex}_3(n, \mathcal{B}(DS_{t,k})),
  \]
   which proves the lower bound.
  For the upper bound, it suffices to deal with the case where $t=k$.
  Like in the proof of Theorem~\ref{thm: spider}, we prove a stronger version for at most $3$-uniform hypergraphs.
  For a $\{2,3\}$-uniform hypergraph $\mathcal{H}$, define $f(\mathcal{H}) = \sum_{e \in E(\mathcal{H})} |e|$.
  Let $\mathcal{H}$ be an at most $3$-uniform counter-example with the minimum value of $f(\mathcal{H})$.
  By a similar argument as in the proof of Theorem~\ref{thm: spider}, we may assume $d^N(u) \ge 2k+1$ for every vertex $u$.
  Now we would like to find a Berge copy of $DS_{k,k}$ in $\mathcal{H}$.
  We start by embedding the edge $uv$ for arbitrary $u$ and $v$.
    
    \begin{claim}\label{claim: for double star}
      Let $\mathcal{H}$ be a $\{2,3\}$-hypergraph.
      Let $u,v$ be two vertices and $\mathcal{E}$ be a set of hyperedges such that $uv$ is not contained in any hyperedge in $\mathcal{E}$.
      If $\left|N_{\mathcal{E}}(u) \right| \ge 2t-1$,  
      $\left|N_{\mathcal{E}}(v) \right| \ge 2s$, then there is a Berge $S_t$ with $u$ being the center vertex and a Berge $S_s$ with $v$ being the center vertex.
      Moreover, the defining vertices are disjoint, the defining hyperedges are in $\mathcal{E}$ and disjoint.
  \end{claim}
  
  \noindent
  \textit{Proof of Claim~\ref{claim: for double star}.}
 We prove this by induction on $t+s$.
  When $t=0$ or $s=0$, the statement can be proved greedily.
  So we may assume $t,s \ge 1$ and the statement is true for all smaller $t+s$.
  For a set of hyperedges $\mathcal{E}$ and a vertex set $W$, define $\mathcal{E}\setminus W = \{e \setminus W \mid e \in \mathcal{E}\}$. When $W$ only contains one vertex $w$, we simply write $\mathcal{E}\setminus w$ for $\mathcal{E}\setminus \{w\}$.
  Let $S_u = N_{\mathcal{E}}(u)$ and $S_v = N_{\mathcal{E}}(v)$.
  
If $S_u \nsubseteq S_v$, pick $w_1 \in S_u \setminus S_v$.
By the definition, there exists a hyperedge $e_{w_1}$ containing $u$ and $w_1$.
  Let $\mathcal{E}' = \mathcal{E}' = \left( \mathcal{E} - \{e_{w_1}\} \right) \setminus w_1$.
  Then $N_{\mathcal{E}'}(u) \supseteq S_u - e_{w_1}$ has size at least $2(t-1)-1$ and $N_{ \mathcal{E}'}(v) = S_v $ has size at least $2s$.
  By induction, we have a Berge $S_{t-1}$ and a Berge $S_{s}$ with $u$ and $v$ as the center vertices respectively.
  Combining with the hyperedge $e_{w_1}$, we have a Berge $S_{t}$ and a Berge $S_s$ with $u$ and $v$ as the center vertices respectively, a contradiction. 
  So $S_u \subseteq S_v$.
  Similarly, we have a contradiction when $S_v \nsubseteq S_u$.

  When $S_u = S_v = S$, then we color the vertex $w$ in $S$ blue if $uw \in \mathcal{E}$, color the pair $ww'$ blue if $uww' \in \mathcal{E}$.
  Similarly, we color the vertex $w$ in $S$ red if $vw \in \mathcal{E}$, color the pair $ww'$ red if $vww' \in \mathcal{E}$.
  Then every vertex in $S$ is either blue~(resp. red) or in a blue~(resp. red) pair.

  If there exists a blue vertex $w$~(that is, $uw \in \mathcal{E}$), then let $\mathcal{E}' = \mathcal{E}\setminus \{w\}$.  
  Then $|N_{\mathcal{E}'}(u)| \ge 2(t-1)$ and $|N_{\mathcal{E}'}(v)| \ge 2s-1$.
  Then we can get a Berge $S_{t-1}$ and a Berge $S_{s}$ using hyperedges in $\mathcal{E}'$ with $u$ and $v$ as the center vertices respectively by induction.
  Combining with the hyperedge $uw$, we have a Berge $S_{t}$ and a Berge $S_s$ with $u$ and $v$ as the center vertices respectively, a contradiction.

  If $w_1w_2$ is a blue pair and $w_1$ is red, then we fix that the hyperedge $uw_1w_2$ corresponds to the edge $uw_2$ and  the hyperedge $vw_1$ corrsponds to the edge $vw_1$.
  Let $\mathcal{E}' = \mathcal{E}\setminus \{w_1,w_2\}$ and then $|N_{\mathcal{E}'}(u)| \ge 2(t-1)-1$ and $|N_{\mathcal{E}'}(v)| \ge 2s-2$.
  By the induction hypothesis, we can get a contradiction.
  Since every vertex is in a blue edge, it actually shows that there is no red vertex.
  Now we have a coloring of edges~(one edge can be both blue and red) in a graph on $S$ satisfying that every vertex is in a blue edge and a red edge.
  
  Assume first that there exists an even cycle such that the edges are alternately blue and red~(a bicolored edge can be viewed as $2$-cycle)~(see Figure~\ref{fig: double star 01}). 
  Assume the even cycle is on vertices $w_1w_2\ldots w_{2\ell}$ where $w_{2p-1}w_{2p}$ is blue and $w_{2p}w_{2p+1}$ is red for $p=1,2,\ldots,\ell$~(let $w_{2\ell+1} = w_1$).
  Then we fix that  the hyperedge $uw_{2p-1}w_{2p}$ corresponds to the edge $uw_{2p-1}$ and the hyperedge $vw_{2p}w_{2p+1}$ corresponds to the edge $vw_{2p}$ with  for $p=1,2,\ldots,\ell$.
  Let $W = \{w_1,w_2,\ldots,w_{2\ell}\}$ and $\mathcal{E}' = \mathcal{E}\setminus W$.
  Then $N_{\mathcal{E}'}(u)\supseteq S\setminus W $ has size at least $2(t-1)-2\ell$ and $N_{\mathcal{E}'}(v)\supseteq S\setminus W $ has size at least $2s-2\ell$.
  Let $t' = \max\{t-\ell, 0\}, s' = \max\{s-\ell, 0\}$.
  By induction, we can get a Berge $S_{t'}$ and a Berge $S_{s'}$ using hyperedges in $\mathcal{E}'$ with $u$ and $v$ as the center vertices respectively.
  Combining $uw_{2p-1}$ and $vw_{2p}$ for $p=1,2,\ldots,\ell$, we have a Berge $S_{t}$ and a Berge $S_s$ with $u$ and $v$ as the center vertices respectively, a contradiction.

  If there is no alternating even cycle, then start with an arbitrary vertex $w_1$, it is contained in a blue edge, say $w_1w_2$.
  Since $w_2$ is also contained in a red edge, say $w_2w_3$, we can continue the process until the alternating path returns to a vertex already on the path.
  Since there is no alternating even cycle, it must result in an ``almost alternating'' odd cycle, i.e., an odd cycle on vertices $w_1w_2\ldots w_{2\ell-1}$ for some $\ell$ such that $w_1w_2, w_3w_4,\ldots,w_{2\ell-1}w_1$ have the same color, say blue, and the rest of the edges are red~(see Figure~\ref{fig: double star 01}).
  Since $w_1$ is contained in a red edge, say $w_1w_1'$, $w_1'$ is not on the cycle, otherwise it would result in an even alternating cycle.
  Similarly to the case of an even alternating cycle, we fix that the hyperedge $vw_1w_1'$ corresponds to the edge $vw_1'$ and the hyperedge $vw_{2p}w_{2p+1}$ corresponds to the edge $vw_{2p}$ with for $p = 1,2,\ldots,\ell-1$.
Finally, the hyperedge $uw_{2p-1}w_{2p}$ corresponds to the edge $uw_{2p-1}$.
  In this case, we also successfully embedded $2\ell$ edges.
  By a similar argument as in the previous case, we can delete these hyperedges and get  a Berge $S_{t'}$ and a Berge $S_{s'}$ with $u$ and $v$ as the center vertices, to obtain a contradiction.
  The details are omitted here.
  \hfill $\blacksquare$ \par

  \begin{figure}[t]
    \centering
    \includegraphics[width=0.8\textwidth]{src/hyperchromatic-15.png}
    \caption{The solid lines and dotted lines represent blue edges and red edges respectively. The left figure shows the case when there is an even cycle with alternating blue and red edges. The right figure shows the case when there is an odd cycle which is almost alternating.}\label{fig: double star 01}
  \end{figure}

  Now we are ready to prove the theorem.
  Recall that $\mathcal{H}$ is a $\{2,3\}$-uniform counter-example with the minimum value of $f(\mathcal{H})$.
  We say a pair is \textit{light} if it is contained in at most one hyperedge, otherwise it is \textit{heavy}.
  Then we claim that every $2$-edge in $\mathcal{H}$ is light, otherwise deleting the $2$-edge does not change the strong chromatic number but results in a counter-example with smaller $f(\mathcal{H})$, a contradiction.
  Similarly, every $3$-edge in $\mathcal{H}$ contains at least two light pairs. Otherwise, if $v_1v_2$ and $v_1v_3$ are heavy in a hyperedge $v_1v_2v_3$, then replacing $v_1v_2v_3$ by $v_2v_3$ does not change the strong chromatic number but results in a counter-example with smaller $f(\mathcal{H})$, a contradiction.

  If there exists a $2$-edge $uv$, then let $\mathcal{E} = E(\mathcal{H}) \setminus \{uv\}$ and then apply Claim~\ref{claim: for double star}.
  We are done in this case.

  Thus, we may assume $\mathcal{H}$ is $3$-uniform. If the hypergraph is not linear, let $uvw$ be a hyperedge such that $uw$ is heavy and $uv, vw$ are light. Then let $\mathcal{E} = E(\mathcal{H}) \setminus \{uvw\}$.
  In this case, $N_{\mathcal{E}}(u)  = N_{\mathcal{H}}(u) \setminus \{v\}$ has size at least $2k$ and $N_{ \mathcal{E}}(v)  \supseteq N_{\mathcal{H}}(v) \setminus \{u, w\}$ has size at least $2k-1$.
  Therefore, we can apply Claim~\ref{claim: for double star} again to get a Berge copy of $DS_{k,k}$.

  Hence, we may assume that $\mathcal{H}$ is $3$-uniform and linear. Then for every vertex $u$, $d^N(u) \ge 2k+1$ actually implies $d^N(u) \ge 2k+2$. Then we pick any $u,v$ contained in a hyperedge $uvw$ and let $\mathcal{E} = E(\mathcal{H}) \setminus \{uvw\}$. Apply Claim~\ref{claim: for double star} again and then we are done.
\end{proof}

\begin{thm}\label{thm: broom}
Let $T$ be a broom $B_{t,k}$~($t,k \ge 2$) which is a tree obtained by identifying an endpoint of $P_t$ with the center vertex of an $S_k$~(see Figure~\ref{fig: special trees}).
$B_{t,k}$ has $t+k$ edges and $t+k+1$ vertices.
When $t \le k-1$
\[
    \mathrm{schex}_3(n, \mathcal{B}(B_{t,k})) = 2k+1.
    \]
    When $t \ge k$, 
\[
    \mathrm{schex}_3(n, \mathcal{B}(B_{t,k})) = t+k.
\]
\end{thm}

\begin{proof}
In the following proof, we always use $v'$ to denote the image of $v$ when we try to embed the $B_{t,k}$.
  When $t \le k-1$, the lower bound comes from $2k+1 = \mathrm{schex}_3(n, \mathcal{B}(S_{k+1})) \le \mathrm{schex}_3(n, \mathcal{B}(B_{t,k}))$ since $S_{k+1}$ is a subgraph of $B_{t,k}$.
 The upper bound can be obtained by a similar proof as Theorem~\ref{thm: main Berge-T, general upper bound}. 
  To be more precise, let us prove the upper bound by proving a stronger statement for at most $\{2,3\}$-uniform hypergraphs as before.
  Let $\mathcal{H}$ be a $\{2,3\}$-uniform counter-example with minimum number of hyperedges, then we may assume that $d^N(u) \ge 2k+1$ for every vertex $u \in V(\mathcal{H})$.
  By Lemma~\ref{lem: minimum degree to berge T}, there is a Berge copy of $S_{k+1}$.
  Assume $v_1'$ is the central vertex of the $S_{k+1}$ with leaves $u_1',u_2',\ldots,u_k',v_2'$.
  Then, since $d^N(v_2) \ge 2k+1$, there exists a vertex $v_3' \notin \{v_1',v_2',u_1',\ldots,u_k'\}$ such that $v_2'v_3'$ is contained in a hyperedge which is different from the defining hyperedges of the $S_{k+1}$.
  We can continue this process and embed the path $v_1'v_2'\ldots v_{t+1}'$ because $t \le k-1$.
  
  When $t \ge k$, the lower bound is given by a clique on $(t+k)$ vertices.
  It remains to prove the upper bound when $t \ge k$.
  Like before, we prove the theorem for $\{2,3\}$-uniform hypergraphs.
  Let $\mathcal{H}$ be a $\{2,3\}$-uniform counter-example with the minimum number of hyperedges.
  We may assume $d^N(u) \ge t+k$ for every vertex $u \in V(\mathcal{H})$.

  First, consider the case when $t \ge 3$.
  We do the embedding process as in Proposition~\ref{prop: berge tree, r uniform} by embedding the star first and then the path.
  The only possible failure is when we embed the last vertex of the path.
  Suppose we have embedded a star $v_1',u_1',\ldots,u_k'$ and a path $v_1',v_2',\ldots,v_t'$ such that for every edge, there exists a distinct hyperedge containing it. Let us denote by $\cB$ this broom.

  Now we would like to find a new vertex $v_{t+1}'$ such that $v_{t}'v_{t+1}'$ is contained in a new hyperedge.
  Since $d^N(v_{t}') \ge t+k$, it only fails when $N_{\mathcal{H}}(v_t') = \{u_1',\ldots,u_k',v_1',\ldots,v_{t-1}',w\}$, where $w$ is the third vertex in the defining hyperedge of $v_{t-1}'v_t'$. 
  In this case, $w$ and $v_t'$ are equivalent; in other words, we could replace $v_{t}'$ by $w$, so the neighbors of $w$ in $\mathcal{H}$ are also exactly $u_1',\ldots,u_k',v_1',\ldots,v_{t-1}',v_t'$.

  Since $w$ is connected to $v_1'$, if the hyperedge containing $v_1'$ and $w$ is the defining hyperedge corresponding to $v_1'u_i'$ for some $i$, or a non-defining hyperedge, then we can replace $u_i'$ by $w$ and let $v_{t+1}' = u_i'$ to get a Berge copy of $B_{t,k}$, a contradiction.
  Therefore, the hyperedge containing $v_1'$ and $w$ is the defining hyperedge corresponding to $v_1'v_2'$. 
  Note that there also exists a hyperedge containing $v_1'$ and $v_t'$, which is not the defining hyperedge of $v_1'v_2'$.
  Since $w$ and $v_t'$ are equivalent, we can exchange the vertices, then the above analysis leads to a contradiction.

    When $t=2$, first there exists a Berge copy of $S_{k+1}$ since $\mathrm{schex}_3(n, \mathcal{B}(S_{k+1}))=2k+1$.
    Assume the $S_{k+1}$ has $v_1'$ as the central vertex and leaves $u_1',u_2',\ldots,u_{k+1}'$.
    For each $u_i'$, since $d^N(u_i') \ge k+2$, there must exist a neighbor of $u_i'$ not in $\{v_1',u_1',\ldots,u_{k+1}'\}$ in $\partial \mathcal{H}$, say $w_i$.
    Then the hyperedge containing $u_i'$ and $w_i$ must be the defining hyperedge of $v_1'u_i'$, otherwise we can extend $u_i'w_i$ to get a $B_{2,k}$.
    Moreover, the neighbors of $u_i'$ in $\partial \mathcal{H}$ are exactly $\{w_i, v_1',u_1',\ldots,u_{k+1}'\}\setminus \{u_i'\}$.
    Then $u_1'u_2'$ is contained in a non-defining hyperedge.
    If we replace $u_2'$ by $w_2$ where the defining hyperedge of $v_1'w_2$ is still $v_1'u_2'w_2$, then we can extend $u_1'u_2'$ with a new hyperedge, which leads to a Berge copy of $B_{2,k}$.
\end{proof}

\section{Strong/Weak chromatic number for Berge-path-free hypergraphs}\label{sec: berge path}

In this section, we study the special case when we forbid Berge copies of a path in $3$-uniform and $4$-uniform cases.
For the $3$-uniform case, we prove Theorem~\ref{thm: main Berge path r=3} to get a sharp bound.
For the $4$-uniform case, we prove Theorem~\ref{thm: main 4uniform berge path} to get a rough bound which is sharp asymptotically.

\subsection{The \texorpdfstring{$3$}{3}-uniform case}\label{subsec: 3uniform path}

When $T = P_k$ is a path with $k$ edges, in the $3$-uniform case, Theorem~\ref{thm: main Berge-T, general upper bound} shows that $\mathrm{schex}_3(n, \mathcal{B}(P_k)) \le k+1$.
To prove Theorem~\ref{thm: main Berge path r=3}, we only need to prove the following proposition which is sufficient.

\begin{prop}\label{prop: strong chromatic number path r=3}
    Let $P_k$ be a path with $k$ edges~($k \ge 3$).
    Let $\mathcal{H}$ be a $\{2,3\}$-hypergraph forbidding Berge copies of $P_k$.
    Then $\chi_s(\mathcal{H}) \le k$.
    Moreover, $\chi_w(\mathcal{H}) \le \lceil k/2 \rceil$.
    The bound is sharp when $\mathcal{H}$ is disjoint union of cliques on $k$ vertices.
\end{prop}

The following lemma is used in many papers studying $\mathcal{B}(P_k)$-free hypergraphs.

\begin{lemma}\label{lemma: berge Ck a component}
  If there exists a Berge copy of $C_k$~(cycle with $k$ edges) in a $\mathcal{B}(P_k)$-free hypergraph $\mathcal{H}$,
  then the defining vertices of the $C_k$ constitute a component of $\mathcal{H}$.
\end{lemma}

\begin{proof}
  Let $U$ be the set of defining vertices of the $C_k$.
  First all the defining hyperedges are contained in $U$, otherwise we have a Berge copy of $P_k$ immediately.
  For a vertex $v \notin U$, if there is a hyperedge containing $v$ and a vertex in $U$, then we have a Berge copy of $P_k$, a contradiction.
\end{proof}

\noindent
\textbf{Proof of Proposition~\ref{prop: strong chromatic number path r=3}.}

  Let $\mathcal{H}$ be a counterexample with minimum number of vertices.
  Then $\mathcal{H}$ must be connected.
  By Lemma~\ref{lemma: berge Ck a component}, we may assume $\mathcal{H}$ is $\mathcal{B}(C_k)$-free.
  Moreover, we may assume $d^N(u) \ge k$ for all $u \in V(\mathcal{H})$.
  Otherwise, we can color $\mathcal{H}\setminus \{u\}$ by induction and then color $u$ greedily.

We prove the statement by induction on $k$.
When $k=3$, it is easy to verify the statement.
From now on, we assume $k \ge 4$ and the statement is true for all $k' < k$.
If there is no Berge path of length $k-1$ in $\mathcal{H}$, we are done by the induction hypothesis.
Otherwise, let $P$ be a Berge path of length $k-1$ on vertices $v_1,v_2,\ldots,v_k$.
Then the neighbors of $v_1$ in $\mathcal{H}$ must be exactly $v_2,v_3,\ldots,v_k$ and $u_1$ where $u_1$ is the third vertex in the defining hyperedge of $v_1v_2$.
Otherwise, if $w$ is a neighbor of $v_1$ in $\mathcal{H}$ outside the path $P$ and $w$ is not in the defining hyperedge of $v_1v_2$, then we have a Berge copy of $P_k$ on $w,v_1,\ldots,v_k$, a contradiction.

Then $v_1$ and $v_k$ are connected by a hyperedge, which must be the defining hyperedge of $v_{k-1}v_k$, otherwise, we have a Berge copy of $C_k$.
Since $v_1$ and $v_k$ are equivalent, the neighbors of $v_k$ must be exactly $v_1,v_2,\ldots,v_{k-1}$ and $u_k$ where $u_k$ is the third vertex in the defining hyperedge of $v_kv_{k-1}$.
Then $u_k = v_1$, a contradiction.
\hfill $\square$ \par

In the graph case, it is known that a $P_k$-free graph has chromatic number at most $k$.
Moreover, there is an efficient algorithm to color a $P_k$-free graph $G$ with $k$ colors.
We can first run depth-first search~(DFS) to generate a DFS-tree.
Since $G$ is $P_k$-free, the DFS-tree has height at most $k-1$.
The we can color the vertices with the same color if they are in the same level.
The time complexity is $O(n+m)$ where $n$ is the number of vertices and $m$ is the number of edges.

Now in the $\{2,3\}$-uniform case, we also propose an efficient algorithm.

\begin{thm}\label{thm: DFS algorithm strong r=3 path}
  Let $P_k$ be a path with $k$ edges~($k \ge 3$).
  Let $\mathcal{H}$ be a $\{2,3\}$-hypergraph forbidding Berge copies of $P_k$ with $n$ vertices and $m$ hyperedges.
  Then there exists an algorithm to find a strong coloring of $\mathcal{H}$ with $k$ colors in $O(n+m)$ time.
\end{thm}

To prove the theorem, we need to do some preparation.
Let $\mathcal{H}$ be a $\mathcal{B}(P_k)$-free $\{2,3\}$-hypergraph.
First we can iteratively remove vertices with degree at most $k-1$ in $\partial \mathcal{H}$ until no such vertex exists.
Let $S$ be the set of deleted vertices and $\mathcal{H}' = \mathcal{H}\setminus S$.
If we can get a strong coloring of $\mathcal{H}'$ with $k$ colors, then we can extend it to a strong coloring of $\mathcal{H}$ with $k$ colors by coloring the vertices in $S$ greedily.
The time complexity is still $O(n+m)$.
So in the following, we may assume every vertex $u$ satisfies $d^N(u) \ge k$.
Moreover, we may assume $\mathcal{H}$ is connected.
Otherwise, we can color each component of $\mathcal{H}$ separately and then combine the colorings without increasing the time complexity.
By Lemma~\ref{lemma: berge Ck a component}, if there is a Berge copy of $C_k$ in $\mathcal{H}$, then $\mathcal{H}$ only contains one component on $k$ vertices, and it is trivial to color.
To summarize, we say $\mathcal{H}$ is \textit{good} if it is connected, $\mathcal{B}C_k$-free and every vertex $u$ satisfies $d^N(u) \ge k$.

  We propose a modified DFS procedure to generate a DFS-tree as Algorithm~\ref{alg: DFS}.
  During the algorithm, for each edge $e$ in the tree, we attach a hyperedge $f(e)$ which contains the pair of vertices in $\mathcal{H}$.
    
    \begin{algorithm}
        \caption{DFS in $\{2,3\}$-hypergraphs}  
        \label{alg: DFS}       
        \begin{algorithmic}[1]    %
          \Require good $\{2,3\}$-hypergraph $\mathcal{H}$, a starting vertex $r$, a set $V$ of visited vertices, a set $\mathcal{E}$ of ignored hyperedges.
          \Ensure A tree $T$ on $V(\mathcal{H})\setminus V$
      
          \Procedure{DFS}{$r$, $V$, $\mathcal{E}$}  
            \State Let $r$ be the root of the tree $T$. Update $V \rightarrow V \cup \{r\}$.
            \While{There exists a vertex $v \in V( \mathcal{H}) \setminus V$ adjacent to $r$ by a hyperedge $e$ not in $\mathcal{E}$}          
              \State Apply $DFS(v, V, \mathcal{E} \cup \{e\})$, get a tree $T'$.
              \State Update $V \rightarrow V \cup V(T')$.
              \State Attach $T'$ to $T$ by connecting $r$ and $v$~(root of $T'$) and let $f(rv) = e$. Let $v$ be the rightmost child of $r$ at this step.
              \If{The third vertex $u$ of $e$~($u = e\setminus\{r,v\}$) is not in $V$}
                \State Apply $DFS(u, V, \mathcal{E})$, get a tree $T''$.
                \State Attach $T''$ to $T$ by connecting $r$ and $u$~(root of $T''$) and let $f(ru) = e$. Let $u$ be the rightmost child of $r$ at this step.
              \EndIf
            \EndWhile
            \State Return $T$.
          \EndProcedure       
        \end{algorithmic}
      \end{algorithm}


      Let $T=DFS(r, \emptyset, \emptyset)$ be the DFS-tree generated by Algorithm~\ref{alg: DFS} on a $\{2,3\}$-hypergraph $\mathcal{H}$ starting from vertex $r$.
      For convenience, we say $T$ is the DFS-tree of $\mathcal{H}$ at $r$.
      From Algorithm~\ref{alg: DFS}, we have the following lemma.


      \begin{lemma}\label{lem: DFS tree special}
        Let $\mathcal{H}$ be a $\{2,3\}$-hypergraph and $T$ be the DFS-tree of $\mathcal{H}$ at $r$.
        If $u_1,u_2$ are two vertices such that the subtrees rooted at $u_1$ and $u_2$ are disjoint, and there exists a hyperedge $e$ containing both $u_1$ and $u_2$,
        then this only happens when $u_1,u_2$ have the same depth $d$ in $T$ with a common parent node $w$. Moreover, the hyperedge $e$ must be exactly $u_1u_2w$. 
        In this case, we say $e$ is a \textbf{special hyperedge} of depth $d$, and $(u_1,u_2)$ is a \textbf{special pair} of depth $d$. 
        Moreover, the special pairs in each depth form a matching.  
        
      \end{lemma}

      \begin{proof}
        Without loss of generality, assume $u_1$ is before $u_2$ in the 
        order of adding vertices to $V$.
        Let $e$ be the hyperedge containing $u_1$ and $u_2$.
        During the algorithm, at some time point, $DFS(u_1, V, \mathcal{E})$ is called for some $V$ and $\mathcal{E}$ to generate the subtree rooted at $u_1$.
        At this time point, $V$ is the set of visited vertices, so $u_2 \notin V$.
        Since $u_1$ and $u_2$ are contained in $e$ and $u_2$ is not in the subtree rooted at $u_1$, it only happens when $e \in \mathcal{E}$.
        During the algorithm, every edge in $\mathcal{E}$ must be $f(v_1v_2)$ for some $v_1,v_2 \in V$.
        Since $e$ contains both $u_1$ and $u_2$ and $\mathcal{H}$ is at most $3$-uniform, it must be $f(u_1u_0)$ where $u_0$ is the parent of $u_1$. 
        By the algorithm, $u_0$ must also be the parent of $u_2$.

        The moreover statement easily follows from the above analysis.
      \end{proof}

      \begin{lemma}\label{lem: DFS tree depth}
        Let $\mathcal{H}$ be a $\{2,3\}$-hypergraph forbidding Berge copies of $P_k$ and $T$ be the DFS-tree of $\mathcal{H}$ at $r$.
        Then the height of $T$ is at most $k-1$.
      \end{lemma}

      \begin{proof}
        Suppose otherwise.
        Let $v_0 = r$, then there exists a path $v_0v_1v_2\ldots v_{k}$ such that $v_i$ has depth $i$ for $i=0,1,2,\ldots,k$ in $T$.
        During the algorithm, for each $v_iv_{i+1}$, there exists a hyperedge $f(v_iv_{i+1})$ which contains $v_i$ and $v_{i+1}$ for $i=0,1,2,\ldots,k-1$.
        By the algorithm, the hyperedges are distinct, then we have a Berge copy of $P_k$ on $v_0,v_1,v_2,\ldots,v_k$, a contradiction.
      \end{proof}
      
      \begin{lemma}\label{lem: DFS tree coloring preparation}
        Let $\mathcal{H}$ be a good $\{2,3\}$-hypergraph without Berge copies of $P_k$ and $T$ be the DFS-tree of $\mathcal{H}$ at $r$.
        Then one of the following statements holds:
        \begin{enumerate}
          \item The height of $T$ is at most $k-2$.
          \item The height of $T$ is $k-1$ but there is no special hyperedge of depth $k-1$.
        \end{enumerate}
      \end{lemma}

      \begin{figure}
        \centering
        \includegraphics[width=0.3\linewidth]{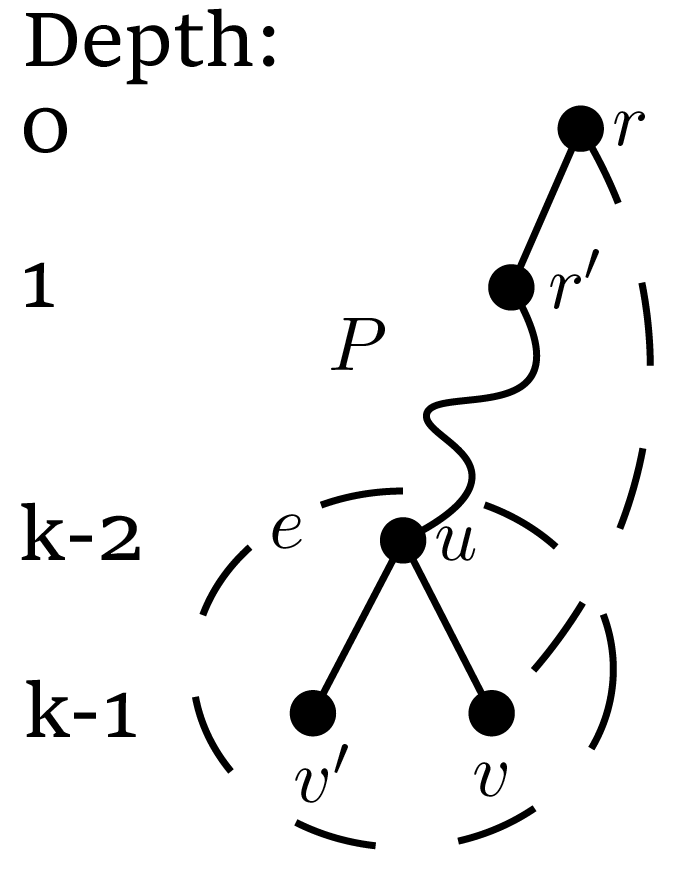}
        \caption{An illustration of the DFS-tree}\label{fig: hyperchromatic 04}
      \end{figure}

      \begin{proof}
        Suppose otherwise.
        Assume $T$ has height $k-1$ and there is a special hyperedge of depth $k-1$, say $v,v'$ with a common parent node $u$.
        Let $P$ be the path from $r$ to $u$ in $T$~(see Figure~\ref{fig: hyperchromatic 04}).
        Note that $P$ has $k-2$ edges and $k-1$ vertices.
        Since for every edge in the path there exists a distinct hyperedge containing its two vertices by the algorithm, this yields a Berge copy of $P_{k-1}$ together with $v$.
        Since $\mathcal{H}$ is good, $v$ has at least $k$ neighbors in $\mathcal{H}$.
        Then the neighbors of $v$ in $\mathcal{H}$ must be exactly the vertices in $P$ and $v'$.
        Since $v$ and $r$ are connected in $\partial \mathcal{H}$, they are contained in a hyperedge $e_{vr}$.
        Since $\mathcal{H}$ is $\mathcal{B}(C_k)$-free, $e_{vr}$ must be the defining hyperedge of $rr'$.

        Since $v$ and $v'$ are equivalent, we can do the same analysis for $v'$ and then obtain a Berge copy of $C_k$, a contradiction.
      \end{proof}

      By Lemma~\ref{lem: DFS tree special}, Lemma~\ref{lem: DFS tree depth} and Lemma~\ref{lem: DFS tree coloring preparation}, we can propose a coloring algorithm based on the DFS-tree as Algorithm~\ref{alg: coloring DFS tree}.
      Note that Lemma~\ref{lem: DFS tree coloring preparation} is an important property of the DFS-tree which guarantees that the coloring algorithm works correctly.

      \begin{algorithm}
        \caption{Coloring the DFS-tree}  
        \label{alg: coloring DFS tree}       
        \begin{algorithmic}[1]    %
          \Require good $\mathcal{B}(P_k)$-free $\{2,3\}$-hypergraph $\mathcal{H}$, the DFS-tree $T$ of $\mathcal{H}$ at $r$, a set $C=\{c_1,c_2,c_3,\ldots,c_t\}$ of $t$ distinct colors.
          \Ensure A colored $T$ with $k$ colors.
      
          \Procedure{Coloring}{$\mathcal{H}$, $T$, $C$}  
            \If {the height of $T$ is at most $2$}
              \State Color $r$ with the color $c_1$.
              \State For every special pair $(u_1,u_2)$ of depth $1$, color $u_1$ and $u_2$ with the color $c_2$ and $c_3$ respectively.
              \State Then color the rest of the vertices with the color $c_2$.
            \Else {~the height of $T$ is at least $3$ and at most $k-1$}
              \State Suppose $r$ is the root of $T$, $u_1,\ldots,u_t,v_1,\ldots,v_t,w_1,\ldots,w_s$ are the children of $r$ in $T$ such that $u_iv_i$ is a special pair for $i=1,2,\ldots,t$ and $w_i$ is not contained in any special pair for $i=1,2,\ldots,s$.
              \State Call $Coloring(\mathcal{H}, T(x), C\setminus \{c_1\})$ for each child $x$ of $r$ where $T(x)$ is the subtree of $T$ rooted at $x$.
              \State Shuffle the colors of each $T(u_i)$ such that $u_i$ and $v_i$ have the different colors for $i=1,2,\ldots,t$.
              \State Color $r$ with the color $c_1$.
            \EndIf
            \State Return the colored $T$.
          \EndProcedure       
        \end{algorithmic}
      \end{algorithm}

    Note that when $T$ has height one, it may require three colors when there is a special pair of depth one.
    If there is no special pair of depth one, then we can color $T$ with two colors.
    So during the algorithm, when $T$ has height $h$, it requires $h+2$ colors when there is a special pair of depth $h$, otherwise, $h+1$ colors are enough.
    Lemma~\ref{lem: DFS tree coloring preparation} guarantees that either $T$ has height at most $k-2$ or there is no special hyperedge of depth $k-1$.
    In either case, the algorithm requires at most $k$ colors.
    It completes the proof of Theorem~\ref{thm: DFS algorithm strong r=3 path}.

    As a side product, the DFS-tree can be used to give a weak coloring as well.
    The situation is easier in this case, we only need to color the vertices in the same depth in $T$ with the same color and every color is used at most twice in two non-consecutive layers.
    It leads to the following theorem.

    \begin{thm}\label{thm: weak coloring DFS algorithm r=3 path}
      Let $\mathcal{H}$ be a good $\mathcal{B}(P_k)$-free $\{2,3\}$-hypergraph.
      Then $\chi_w(\mathcal{H}) \le \lceil k/2 \rceil$.
      Moreover, there exists an algorithm to color $\mathcal{H}$ with at most $\lceil k/2 \rceil$ colors in time $O(n+m)$.
    \end{thm}

\subsection{The \texorpdfstring{$4$}{4}-uniform case}\label{subsec: 4uniform path}

In this subsection, we prove Theorem~\ref{thm: main 4uniform berge path} by proving Proposition~\ref{prop: berge path 4 uniform}, which is sufficient.

\begin{prop}\label{prop: berge path 4 uniform}
  Let $\mathcal{H}$ be an at most $4$-uniform hypergraph without Berge copies of $P_k$.
  Then $\chi_s(\mathcal{H}) \le k + 9$.
  It implies that $\mathrm{schex}_4(n, \mathcal{B}(P_k)) \le k + 9$.
\end{prop}

\begin{proof}
  We prove the proposition by induction on $k$.
  When $k\le 2$, the proposition is trivial.
  From now on, we may assume $k \ge 3$.

  Then we prove the proposition for a fixed $k$ by induction on $n$.
  When $n \le k+9$, the proposition is trivial.
  Suppose $n > k+9$.
  Let $\mathcal{H}$ be a counterexample with the minimum number of vertices.
  For any vertex $u$, if $d^N(u) < k+9$, then let $\mathcal{H}' = \mathcal{H}\setminus u$.
  By induction, $\chi_s(\mathcal{H}') \le k + 9$.
  Then the coloring of $\mathcal{H}'$ can be extended to $\mathcal{H}$ by greedily coloring $u$ with a color that does not appear in $N_{\mathcal{H}}(u)$, a contradiction.
  So we have $d^N(u) \ge k+9$ for every vertex $u$.

  Since $\mathcal{H}$ is a minimum counterexample, it must be connected.
  By Lemma~\ref{lemma: berge Ck a component}, we may assume $\mathcal{H}$ is $\mathcal{B}(C_k)$-free.
  We further claim that there is no Berge copy of $C_{k-1}$.

  \begin{claim}\label{claim: berge Ck-1}
    There is no Berge copy of $C_{k-1}$ in $\mathcal{H}$.
  \end{claim}

  \begin{figure}[t]
    \centering
    \includegraphics[width=0.7\textwidth]{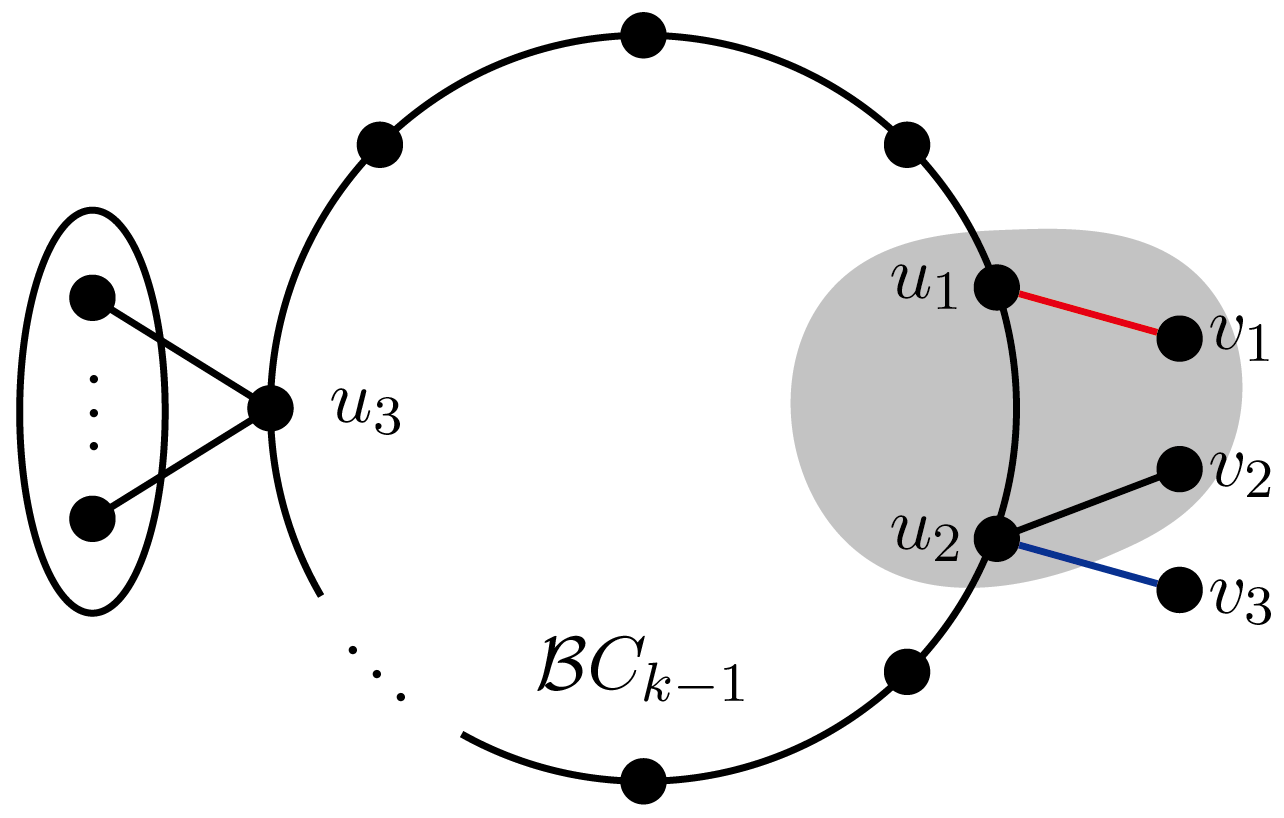}
    \caption{Two good vertices.}\label{fig: berge Ck-1 01}
  \end{figure}

  \noindent
  \textit{Proof of Claim~\ref{claim: berge Ck-1}}
  Suppose there exists a Berge copy of $C_{k-1}$ in $\mathcal{H}$ on vertices $w_1,w_2,\ldots,w_{k-1}$.
  We say a vertex $w_i$ is \textit{good} if there exists at least three vertices $v_1,v_2,v_3$ outside the cycle such that each $w_iv_j$ is contained in a non-defining hyperedge.
  Otherwise, we say it is \textit{bad}. 

  We claim that there are no two consecutive good vertices on the cycle.
  Suppose otherwise.
  Assume $u_1$ and $u_2$ are two consecutive good vertices~(see Figure~\ref{fig: berge Ck-1 01}).
  Then there exists a vertex $v_1$ outside the cycle such that $v_1u_1$ is contained in a non-defining hyperedge.
  Then since there are at least three vertices outside the cycle contained in a non-defining hyperedge together with $u_2$, there must exist a vertex $v_3$ outside the cycle such that $v_3u_2$ is contained in another non-defining hyperedge.
  Then we have a Berge copy of $P_k$, a contradiction.

  Let $\mathcal{P}$ be the set of pairs $(w_i, v)$ where $w_i$ is on the cycle and $v$ is outside the cycle such that there is a defining hyperedge containing $w_i$ and $v$. 
  Then $|\mathcal{P}| \le 4(k-1)$ since $\mathcal{H}$ is at most $4$-uniform.
  Let $u_3$ be a bad vertex on the cycle~(see Figure~\ref{fig: berge Ck-1 01}).
  Note that $d^N(u_3) \ge k+9$, thus $u_3$ has at least 11 neighbors not on the cycle. At most two of those neighbors are in a non-defining hyperedge toether with $u_3$, thus there are at least $9$ pairs $(u_3, v)$ in $\mathcal{P}$ for some $v$ outside the cycle.
  Since there are no two consecutive good vertices, there are at least $\left \lceil \frac{k-1}{2}  \right \rceil$ non-good vertices.
  Then we have $9 \left \lceil \frac{k-1}{2}  \right \rceil   > 4(k-1)$ pairs in $\mathcal{P}$, a contradiction.
  \hfill $\blacksquare$ \par

  By induction on $k$, we may assume that the longest Berge path in $\mathcal{H}$ has length $k-1$ on the vertices $v_1,v_2,\ldots,v_k$.
  Then $v_1$ has at least $10$ neighbors in $\mathcal{H}$ outside the path.
  All such neighbors can be contained only in the defining hyperedges of the path together with $v_1$, otherwise, we have a Berge copy of $P_{k}$, a contradiction.
  Now we define functions $\ell(v_i)$ and $r(v_i)$ as follows~(see Figure~\ref{fig: berge Pk 02}):
  \begin{enumerate}
    \item If $v_i$~($i>1$) is a neighbor of $v_1$ and $v_1,v_i$ are contained in a non-defining hyperedge, then let $\ell_1(v_i) = 1$.
    \item If $v_i$~($i<k$) is a neighbor of $v_k$ and $v_k,v_i$ are contained in a defining hyperedge, then let $r_1(v_i) = 1$.
    \item If the defining hyperedge of $v_iv_{i+1}$ contains both $v_1$ and $v_k$, then let $\ell_2(v_i) = \ell_2(v_{i+1}) = 1/2$, $r_2(v_i) = r_2(v_{i+1}) = 1/2$.
    \item If the defining hyperedge of $v_iv_{i+1}$ contains only $v_1$ but not $v_k$, then let $\ell_3(v_i) =1/2,  \ell_3(v_{i+1}) = 1$.
    \item If the defining hyperedge of $v_iv_{i+1}$ contains only $v_k$ but not $v_1$, then let $r_3(v_i) = 1,  r_3(v_{i+1}) = 1/2$.
    \item Finally, let $\ell(v_i)$ be the maximum amount of $\ell_1(v_i), \ell_2(v_i), \ell_3(v_i)$ and $r(v_i)$ be the maximum amount of $r_1(v_i), r_2(v_i), r_3(v_i)$. For convenience, let $\ell(u) = r(u) = 0$ if $u$ is not on the path and let $\ell(v_1) = r(v_k) = 0$.
  \end{enumerate}

  \begin{figure}[t]
    \centering
    \includegraphics[width=0.7\textwidth]{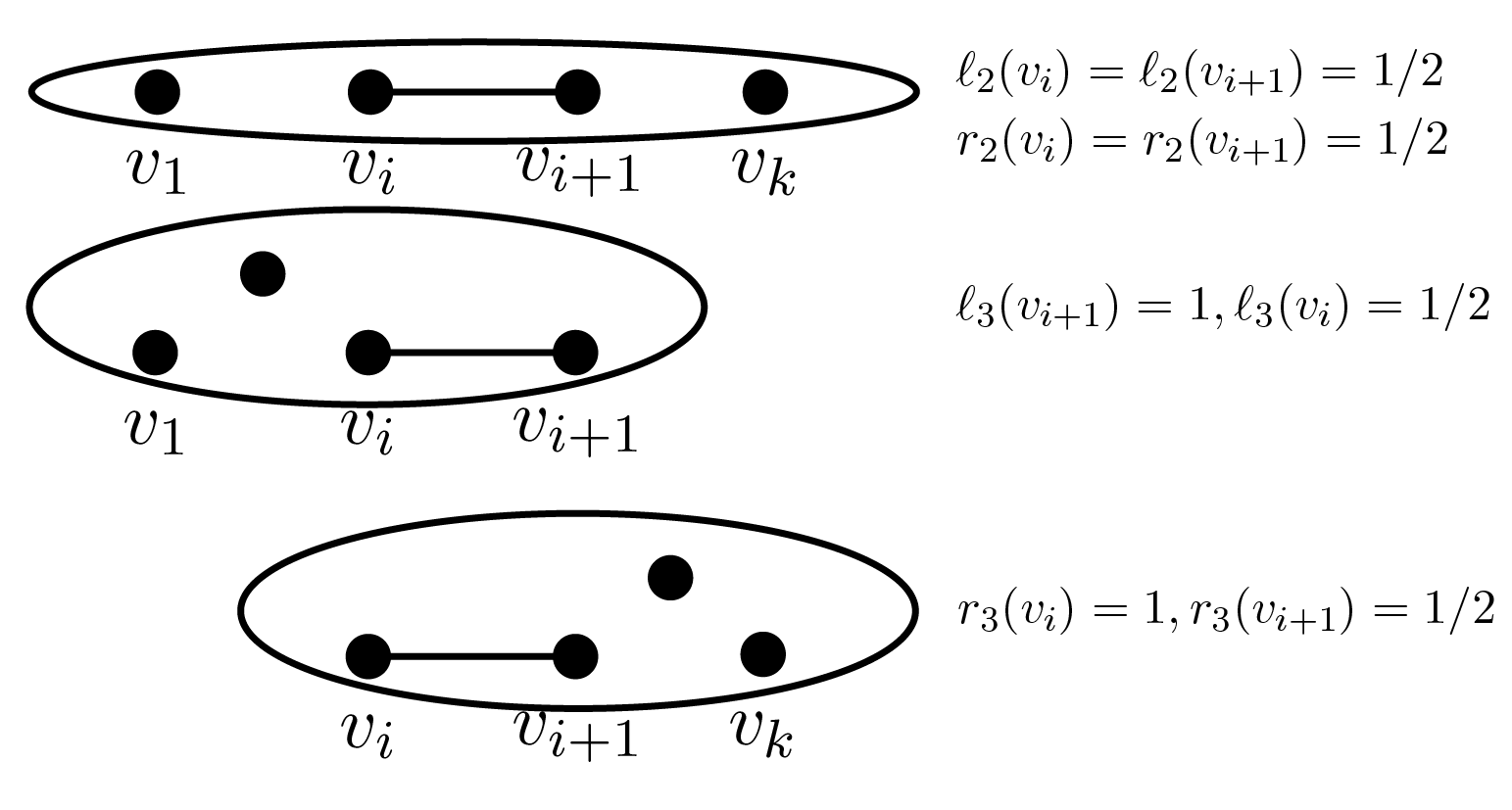}
    \caption{The function $\ell(v_i)$ and $r(v_i)$ for $v_i$.}\label{fig: berge Pk 02}
  \end{figure}

  \begin{claim}\label{claim: half weight}
    We claim that $\sum_{i=2}^{k} \ell(v_i) \ge k/2+4$ and $\sum_{i=1}^{k-1} r(v_i) \ge k/2+4$.
  \end{claim}

  \noindent
  \textit{Proof of Claim~\ref{claim: half weight}.}
  We only prove the first inequality, the second inequality can be proved similarly.
  Let $S$ be the set of neighbors of $v_1$ except the vertices in the defining hyperedge of $v_1v_2$, then $|S| \ge k+8$.
  There are only three possible values for $\ell(s)$: $0, 1/2, 1$.
  We claim that a vertex $s \in S$ has $\ell(s) = 0$ only when $s$ is the fourth vertex in the defining hyperedge of $v_{i}v_{i+1}$ for some $1<i<k$. 
  Indeed, for a vertex $s$ with $\ell(s)=0$, it is contained in a hyperedge with $v_1$.
  Since $v_1v_2\ldots v_k$ form a longest Berge-path and by the definition of $S$, the hyperedge must be the defining hyperedge of $v_iv_{i+1}$ for some $1 < i < k$.
  Moreover, if $s \in \{v_i,v_{i+1}\}$, then $\ell(s) \neq 0$ by the rules, so $s$ must be the fourth vertex in the defining hyperedge.

  In this case, let $f(s) = v_{i+1}$, then we have that $f$ is an injective function and $\ell(f(s)) = 1$.
  So in average, we have $\sum_{s \in S} \ell(s) \ge |S|/2 \ge (k+8)/2$.
  \hfill $\blacksquare$ \par

  \begin{claim}\label{claim: weight > 1}
    There exists $v_i$ and $v_{i+1}$ such that $\ell(v_{i+1}) + r(v_{i}) > 1$.
  \end{claim}

  \noindent
  \textit{Proof of Claim~\ref{claim: weight > 1}.}
  Let $S_i$ be the set of vertices $v_j$ such that $\ell(v_j) = i/2$ for $i=0,1,2$.
  Then by Claim~\ref{claim: half weight}, we have $ |S_1|/2 + |S_2| \ge (k+8)/2$. 
  Suppose the statement does not hold, then for each $v_j \in S_i$~($j>1$), $r(v_{j-1}) \le (2-i)/2$.
  Then we have $\sum_{j=1}^{k} r(v_{j}) \le |S_0|-1 + |S_1|/2$~(here the minus one corresponds to $v_1 \in S_0$).
  Combining that with $|S_0| +|S_1| + |S_2| = k$, we have that $\sum_{j=1}^{k} r(v_{j}) \le k/2+3$, a contradiction.
  \hfill $\blacksquare$ \par

Let $v_i$ and $v_{i+1}$ be vertices given by Claim~\ref{claim: weight > 1}.
  Without loss of generality, assume that $\ell(v_{i+1}) = 1$ and $r(v_i) \ge 1/2$.
  Since $\mathcal{H}$ is $\mathcal{B}(C_k)$-free and $\mathcal{B}(C_{k-1})$-free, $v_{i+1}$ cannot be $v_k$.
  By the definition of $\ell(v_i)$, there are only two possible cases when $\ell(v_{i+1}) = 1$:

  \textbf{Case 1:} $v_{i+1}$ and $v_1$ are contained in a non-defining hyperedge. 
  If $v_i$ and $v_k$ are contained in another non-defining hyperedge, then we have a Berge copy of $C_k$: $v_1,\ldots,v_i,v_k,\ldots,v_{i+1}, v_1$.
  If $v_i$ and $v_k$ are contained in the same non-defining hyperedge as $v_1$ and $v_{i+1}$, then 
  this hyperedge extends our Berge path to a Berge copy of $C_k$, a contradiction.
  If $v_i$ and $v_k$ are contained in the defining hyperedge corresponding to $v_{i}v_{i+1}$, then we have a Berge copy of $C_{k}$: we go along our path from $v_1$ to $v_i$, use the defining hyperedge corresponding to $v_{i}v_{i+1}$ to go to $v_k$, then go back along our path to $v_{i+1}$, and then use the non-defining hyperedge to go to $v_1$.
  If $v_i$ and $v_k$ are contained in the defining hyperedge of $v_{i-1}v_{i}$, then we have a Berge copy of $C_{k-1}$ doing the same but skipping $v_i$.
  In all these subcases, we have a Berge copy of $C_k$ or $C_{k-1}$, while in other cases $r(v_i) =0$, a contradiction.

  \begin{figure}
    \centering
    \includegraphics[width=0.7\textwidth]{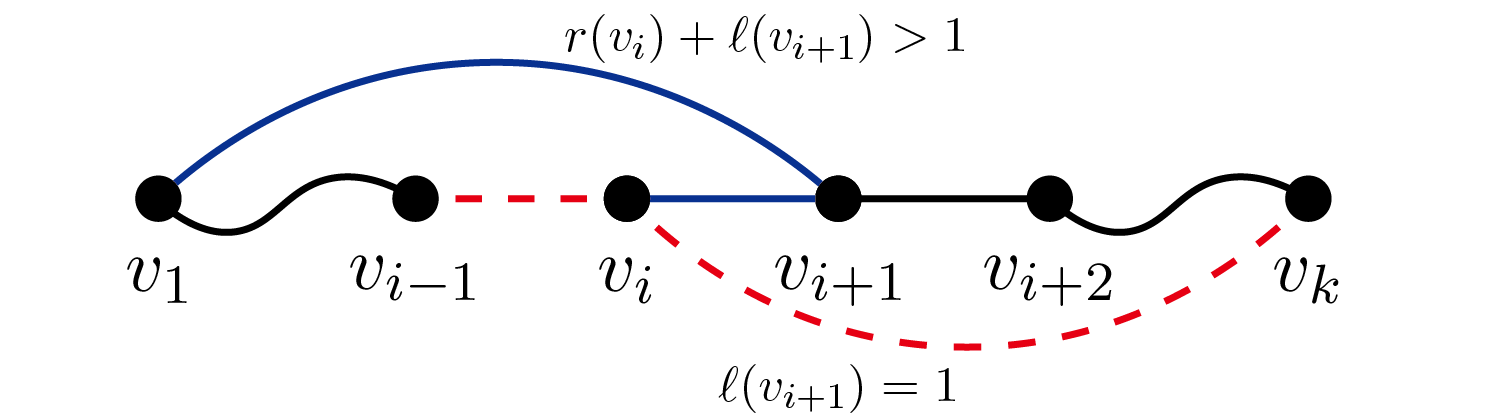}
    \caption{The blue edges~(solid lines) are the defining hyperedge of $v_iv_{i+1}$ and the red edges~(dotted lines) are the defining hyperedge of $v_{i-1}v_{i}$.}\label{fig: berge Pk 03}
  \end{figure}

  \textbf{Case 2:} $v_1$ is contained in the defining hyperedge corresponding to $v_{i}v_{i+1}$ which does not contain $v_k$.
  If $v_i$ and $v_k$ are contained in a non-defining hyperedge, then we have a Berge copy of $C_k$: We go from $v_1$ to $v_i$ along our path, then use this non-defining hyperedge to go to $v_k$, go back to $v_{i+1}$ along our path, and then we go to $v_1$ using the defining hyperedge corresponding to $v_{i}v_{i+1}$.
  Then $r(v_i) \ge 1/2$ only happens when $v_k$ is contained in the defining hyperedge that corresponds to $v_{i-1}v_{i}$~(see Figure~\ref{fig: berge Pk 03}).
  Then we have a Berge copy of $C_{k-1}$ doing the same but skipping $v_i$.

  In conclusion, we can get a contradiction in all cases. This completes the proof.
\end{proof}

\bigskip
\textbf{Funding}: The research of Gerbner is supported by the National Research, Development and Innovation Office---NKFIH under the grant KKP-133819 and by the J\'anos Bolyai scholarship. 

The research of Wang is supported by the China Scholarship Council (Grant 202506210200) and the National Natural Science Foundation of China (Grant 12571372).

The research of Duan is supported by the National Natural Science Foundation of China
(Nos. 12131013, 12171393, 12471334) and Shaanxi Fundamental Science Research Project for Mathematics and Physics (No. 22JSZ009).

\end{document}